\documentclass[amssymb,11pt]{article}
\usepackage{amsfonts}
\usepackage{amsthm}
\usepackage{amsmath}
\usepackage{graphicx}
\usepackage[psamsfonts]{amssymb}
\usepackage{amscd}

\title{$C^1$-generic conservative diffeomorphisms have trivial centralizer}
\author{C. Bonatti, S. Crovisier and A. Wilkinson}

 \def\NN{{\mathbb N}}  
 \def\RR{{\mathbb R}}  
   
 \def\ZZ{{\mathbb Z}}

    \def\cS{{\cal S}}
\def\cB{{\cal B}}    
\def\cC{{\cal C}}   \def\cO{{\cal O}} \def\cU{{\cal U}}
\def\cD{{\cal D}}   \def\cP{{\cal P}} 
\def\cE{{\cal E}}    
   \def\cR{{\cal R}}

\newtheorem{mainthm}{Theorem}

\newtheorem{ques}{Question}

\newtheorem{theorem}{Theorem}[section]
\newtheorem*{thmempty}{Theorem}

\newtheorem{prop}[theorem]{Proposition}

\newtheorem{lemma}[theorem]{Lemma}

\newtheorem{defi}[theorem]{Definition}

\newtheorem{clai}{Claim}
\newtheorem{affi}{Claim}

\theoremstyle{remark}
\newtheorem{remark}[theorem]{Remark}
\newtheorem{remarks}[theorem]{Remarks}

\def\dim{\hbox{dim} }
\def\interior{\hbox{Int} }
\def\sign{\hbox{Sign} }

\def\Lip{\hbox{Lip} }

\def\Diff{\hbox{Diff} }

\def\title{\em}
\def\eps{\varepsilon}

\def\Jac{\hbox{Jac}}

\def\transverse{\,\raise2pt\hbox to1em{\hfil$\top$\hfil}\hskip -1em \hbox
to1em{\hfil$\cap$\hfil}\,}

\newlength{\figboxwidth} 

\begin{document}

\maketitle
\begin{abstract}
We prove that the spaces of $C^1$ symplectomorphisms
and of $C^1$ volume-preserving diffeomorphisms of connected manifolds
both contain residual subsets of diffeomorphisms
whose centralizers are trivial.
\smallskip

\smallskip

\begin{center}
\bf Les diff\'eomorphismes conservatifs $C^1$-g\'en\'eriques ont un centralisateur trivial
\end{center}

\noindent{\it R\'esum\'e~:} Nous montrons que l'espace des symplectomorphismes de classe $C^1$
et l'espace des diff\'emomorphismes de classe $C^1$ pr\'eservant une forme volume contiennent tous
deux des sous-ensembles r\'esiduels de diff\'eomorphismes dont le centralisateur est trivial.
\smallskip

\noindent{\it Mots cl\'e~:} Centralisateur trivial, sym\'etries triviales, propri\'et\'es $C^1$-g\'en\'eriques.
\end{abstract}

\section*{Introduction}

Let $M$ be a connected compact manifold.
The {\em centralizer} of a $C^r$ diffeomorphism $f\in \Diff^r(M)$ is defined 
as $$Z^r(f):=\{g\in \Diff^r(M): fg=gf\}.$$
Clearly $Z^r(f)$ always contains the group $<f>$ of all the powers of $f$.  We 
say that $f$ has {\em trivial centralizer} if $Z^r(f) = <f>$.
A diffeomorphism $f$ with trivial centralizer posesses no smooth symmetries,
such as those that would arise if, for example, $f$ embedded in a flow or
were the lift of another diffeomorphism. S. Smale asked the following:

\begin{ques}[\cite{Sm1,Sm2}]\label{q.smale}
Consider the set of $C^r$ diffeomorphisms of a compact connected manifold $M$ with trivial centralizer.
\begin{enumerate}
\item Is this set dense in $\Diff^r(M)$?
\item Is it {\em residual} in $\Diff^r(M)$?  That is, does it contain a dense $G_\delta$ subset?
\item Does it contain an open and dense subset of $\Diff^r(M)$?
\end{enumerate}
\end{ques}

We refer to~\cite{BCW3} for a discussion on this question.
To summarize,  we mention some cases in which it has been answered:
\begin{itemize} 
\item N. Kopell \cite{Ko} solved the one-dimensional smooth case:
the set of diffeomorphisms with trivial centralizer contains
an open and dense subset of $\Diff^r(S^1)$.
\item For $r\geq 2$ on higher dimensional manifolds, there exists some results
about the set of diffeomorphisms with trivial centralizer under additional
dynamical assumptions. For instance, this set contains
an open and dense subset among the set of Axiom A diffeomorphisms in $\Diff^\infty(M)$ 
possessing at least one periodic sink or source \cite{PY1} or defined on a surface \cite{F}, 
and generic among those satisfying the no-cycles condition \cite{F};
and it is locally residual among a class of partially hyperbolic $C^\infty$ diffeomorphisms
with $1$-dimensional center~\cite{Bu1}.
\item In the $C^1$ setting, Y. Togawa proved \cite{To1,To2} that generic Axiom A diffeomorphisms
have trivial centralizer.
\end{itemize}

This paper deals with the first two parts of Question \ref{q.smale} for all compact $M$
in the case of volume-preserving and symplectic $C^1$ diffeomorphisms.
If $M$ is a symplectic manifold, then $\hbox{Symp}^1(M)$ denotes the space of
$C^1$ symplectomorphisms of $M$.
If $M$ carries a volume $\mu$, then we denote by $\Diff^1_{\mu}(M)$
the space of $C^1$ diffeomorphisms of $M$ that preserve $\mu$.
The spaces $\Diff^1(M)$, $\hbox{Symp}^1(M)$, $\Diff^1_{\mu}(M)$ are Baire spaces in the $C^1$
topology. Recall that a {\em residual} subset of a Baire space is one
that contains a countable intersection of open-dense sets.

\begin{mainthm}\label{t=main1}
Let $M$ be a compact, connected manifold of dimension at least~$2$.
Then:
\begin{enumerate}
\item[(a)] In ${Symp}^1(M)$ the set of diffeomorphisms $f$ whose centralizer $Z^1(f)\cap {Symp}^1(M)$
is trivial is a residual subset.
\item[(b)] In $\Diff^1_{\mu}(M)$ the set of diffeomorphisms $f$ whose centralizer $Z^1(f)\cap \Diff^1_{\mu}(M)$
is trivial is a residual subset.
\end{enumerate}
\end{mainthm}

In a recent work~\cite{BCW2} we also solved Questions~\ref{q.smale} a) and b) in the $C^1$ dissipative case:
the set of diffeomorphisms $f$ whose centralizer $Z^1(f)$ is trivial is a residual subset of $\Diff^1(M)$.
On the other hand with G. Vago~\cite{BCVW} we answered negatively to Questions~\ref{q.smale} c)
in the dissipative and in the symplectic case:
for any compact symplectic manifold $M$,
there exists a familly of symplectomorphisms of $M$ with large centraliser $Z^1(f)\cap {Symp}^1(M)$
and that is dense in a non-empty open subset $U$ of ${Symp}^1(M)$; with a different method we also constructed 
a family of diffeomorphisms of $M$ with large centraliser $Z^1(f)$
that is dense in a non-empty open subset $U$ of $\Diff^1(M)$.
The existence of such examples in the volume-preserving setting in dimension larger or equal to $3$ remains open.

In~\cite{Ko,PY1,PY2,Bu1,To1,To2}, the proof of the triviality of the centralizer
splits in two parts: the first implies that,
for a generic $f$ and for any diffeomorphism $g$ that commutes with $f$,
we have $g=f^ \alpha$ on a dense subset of $M$, where $\alpha$ is a function
that is $f$-invariant; the second part considers the global dynamics
and proves that $\alpha$ is bounded, then constant. It is for this second part that
additional assumptions such as hyperbolicity are often required.
Let us mention that for conservative diffeomorphisms of the $2$-sphere,
L. Burslem has shown~\cite{Bu2} a partial result for the topology $C^r$ with $r\geq 16$
(more precizely, she obtained the first step).

The proof of theorem~\ref{t=main1} also decomposes in two steps. The first one is
provided by the following result.

\begin{mainthm}\label{t=main} (a) For any diffeomorphism $f$ in a residual set $\cR\subset \Diff^1(M)$,
any periodic point $x\in Per(f)$ is hyperbolic and for any $g\in Z^ 1(f)$ there exists $n \in\ZZ$
such that $g$ coincides with $f^n$ on $W^s(x)$.

\noindent (b) For any diffeomorphism $f$ in a residual set $\cR_{\hbox{symp}}\subset \hbox{Symp}^1(M)$,
any hyperbolic periodic point $x\in Per(f)$ and any $g\in Z^ 1(f)$, there exists $n\in\ZZ$
such that $g$ coincides with $f^n$ on $W^s(x)$.

\noindent(c) For any diffeomorphism $f$ in a residual set $\cR_\mu\subset \Diff^1_{\mu}(M)$,
any hyperbolic periodic point $x\in Per(f)$ and any $g\in Z^1(f)$, there exists $n\in\ZZ$
such that $g$ coincides with $f^n$ on either $W^s(x)$ or $W^u(x)$.
\end{mainthm}

The second step is provided by the following genericity result about conservative diffeomorphisms.
This immediately concludes the proof of theorem~\ref{t=main1}.

\begin{thmempty}[\cite{BC,ABC}] For any compact connected manifold $M$, there are residual sets 
$\widetilde\cR_{\hbox{symp}}\subset \hbox{Symp}^1(M)$ and $\widetilde\cR_\mu\subset \Diff^1_\mu(M)$
such that,  every $f\in \widetilde\cR_{\hbox{symp}}\cup \widetilde\cR_\mu$
has a hyperbolic periodic point $p$ with
$$\overline{W^s(p)} = \overline{W^u(p)} = M.$$ 
\end{thmempty}

More generally, Theorem~\ref{t=main} naturally applies to the class of 
$C^1$ diffeomorphisms satisfying a property we call {\em periodic 
accessibility}: this property says that any two points $x,y$ in a 
dense subset $\cE\subset M$ of non-periodic points
may be joined by a finite sequence $x_0=x,x_1,\dots,x_n=y$, $x_i\in \cE$ such that
for any $i\in\{0,\dots,n-1\}$, the points $x_i$ and $x_{i+1}$ belong to the closure
$\overline{W^s(\cO_i)}$ of the stable manifold or to the closure $\overline{W^u(\cO_i)}$ of the
unstable manifold of a hyperbolic periodic orbit $\cO_i$.

\begin{ques} Is periodic accessibility generic in $\Diff^r(M)$?
\end{ques}
As a weaker problem\footnote{This problem is often proposed as a conjecture by the first author.}, one can also ask if, for generic diffeomorphims,
the union of the stable manifolds of the periodic points are dense in $M$.
\bigskip

As an immediate corollary of Theorem~\ref{t=main},
it follows that if $f\in \cR$ satisfies the periodic accessibility property, then $Z^1(f)$ is trivial.
Since the periodic accessibility property is satisfied by Axiom A
diffeomorphisms, one can recover in this way Togawa's result
that the $C^1$ generic Axiom A diffeomorphism has trivial centralizer.
The proof in the general case~\cite{BCW3} is much more delicate and involves new ideas.

Theorem~\ref{t=main}~(a) was previously proved by Togawa \cite{To1,To2}.
We present here a variation of his argument, based on an unbounded distortion property
that is discussed in section~\ref{s=distortion}; we note below that this proof also handles 
the ``Lipschitz centralizer" and this will be crucial for~\cite{BCW3}.
When we started with this work we were not aware of Togawa's papers and
we used a stronger unbounded distortion property which is much more difficult to obtain.
We think that our former approach has independent interest and could have further application:
it is written in the first version~\cite{BCW1} of this article.
The conservative cases (b) and (c) of Theorem~\ref{t=main} can be derived as for case (a)
with an extension result presented in section~\ref{s.invmanif}:
any perturbation of the dynamics of a symplectomorphism inside the stable manifold of a periodic point
can be realized as a perturbation of the dynamics on $M$.

\section{The unbounded distortion property}\label{s=distortion}
Kopell's proof in \cite{Ko} that the set of diffeomorphisms having a trivial centralizer
is open and dense in $\Diff^r(S^1)$ for $r\geq 2$  uses the fact that
a $C^2$ diffeomorphism $f$ of $[0,1]$  without fixed points in $(0,1)$
has {\em bounded distortion}, meaning: for any $x,y\in (0,1)$, the ratio
\begin{eqnarray}\label{e=xydistort}
\frac{|{f^n}'(x)|}{|{f^n}'(y)|}
\end{eqnarray}
is bounded, independent of $n$ and uniformly for $x,y$ lying in a
compact set.  A bounded distortion estimate 
lies behind many results about $C^2$, hyperbolic 
diffeomorphisms of the circle and 
codimension-$1$ foliations.

Suppose that $r\geq 2$. 
Since Morse-Smale diffeomorphisms are open and dense in $\Diff^r(S^1)$,
the proof that the set of diffeomorphisms having a trivial centralizer is open and dense in $\Diff^r(S^1)$ 
essentially reduces to showing that ($C^r$ open and densely)
a $C^r$ diffeomorphism $f:[0,1]\to [0,1]$  without fixed points in $(0,1)$
has trivial centralizer.  The 
bounded distortion of such an $f$ forces its centralizer to
embed simultaneously in two smooth flows containing $f$, one determined
by the germ of $f$ at $0$, and the other by the germ at $1$; 
for an open and dense set of $f\in \Diff_+^r[0,1]$, these flows
agree only at the iterates of $f$.
The $r\geq 2$ hypothesis is clearly necessary for bounded distortion.

The central observation and starting point of this paper is that
the centralizer of a $C^1$ diffeomorphism of $[0,1]$ with 
{\em unbounded} distortion is always trivial.  We elaborate a bit on this.
Notice that if $x$ and $y$ lie on the same $f$-orbit, then the 
ratio in (\ref{e=xydistort}) is bounded, independent of $n$. 
We show that, $C^1$ generically among the diffeomorphisms of $[0,1]$ without
fixed points $(0,1)$, the ratio (\ref{e=xydistort}) is uniformly bounded in $n$ {\em only if}
$x$ and $y$ lie on the same orbit; that is, for a residual set of $f$,
and for all $x,y\in (0,1)$, 
if $x\notin {\mathcal O}_f(y) = \{f^n(y)\,\vert\, n\in \ZZ\}$,
then
\begin{eqnarray}\label{e=xylimsup}
\limsup_{n\to\infty} \frac{|{f^n}'(x)|}{|{f^n}'(y)|} = \infty.
\end{eqnarray}

 Assume that this unbounded  distortion property holds for $f$.  
Fix $x\in(0,1)$.  A simple application of the Chain
Rule shows that if $gf = fg$, then the distortion in (\ref{e=xydistort}) 
between  $x$ and $y=g(x)$ is bounded; hence $x$ and $g(x)$ must
lie on the same $f$-orbit. From here, it is straightforward to show that 
$g= f^n$, for some $n$ (see Lemma~\ref{l.local} below).  
As in \cite{Ko}, a small amount of additional work shows that a residual set in $\Diff^1(S^1)$
has trivial centralizer.  

The bulk of this section is devoted to formulating and 
proving a higher-dimensional version of the argument we have just 
described. The interval is replaced by an invariant manifold (stable or 
unstable) of a periodic point.  The derivative $f'$ in 
(\ref{e=xylimsup}) is replaced by
the Jacobian of $f$ along the invariant manifold.

\subsection{Unbounded distortion along stable manifolds of hyperbolic periodic points}

Denote by $\Lip(M)$ is the set of {\em lipeomorphisms} of $M$:
these are the homeomorphisms $g$ of $M$ such that $g$ and $g^{-1}$ are Lipschitz.
For any hyperbolic periodic orbit $\cO\subset M$ of a diffeomorphism $f\in \Diff^1(M)$
and for $x\in W^s(\cO)$ we denote by $\Jac^s(f)(x)$
the Jacobian of the map induced by $T_xf$ between $T_xW^s(\cO)$ and
$T_{f(x)}W^s(\cO)$. 

\begin{defi}\label{d.distortion}
A hyperbolic periodic orbit $\cO\subset M$ of a diffeomorphism $f\in \Diff^1(M)$
has the {\em stable manifold unbounded distortion property}
if there exists a dense set $\cD\subset W^s(\cO)$ with the following property:
for any points $x\in \cD$ and $y\in W^s(\cO)\setminus \cO$ not in the same $f$-orbit
and for any $N\geq 1$ there exists $n\geq 1$ such that
$$\left|\log(\Jac^s(f^n)(x))-\log(\Jac^s(f^n)(y))\right| = \infty.$$
\end{defi}

As mentioned in the previous subsection, unbounded distortion forces trivial 
centralizers along the stable manifold:

\begin{lemma}\label{l.local}
Let $f\in \Diff^1(M)$ be a diffeomorphism and
$\cO$ be a hyperbolic periodic orbit having the stable manifold unbounded distortion property.
Then for any $g\in \Lip(M)$ such that $g\circ f=f\circ g$ and $g(\cO) = \cO$,
and for any connected component $W$ of $W^s(\cO)\setminus \cO$,
there exists an $i\in \ZZ $ such that $g = f^{i}$ on $W$.
\end{lemma}
\begin{remark}
Since $\cO$ is hyperbolic, if $g$ is differentiable at $p$, then
there exists $i\in \ZZ$ such that $g=f^i$ on $W^s(\cO)$.
\end{remark}

\begin{proof}
Let $\cD\subset W^s(\cO)$ be the dense set as in definition~\ref{d.distortion}.

We claim that:
\emph{For every $x\in \cD\setminus \cO$, the map $g$ preserves the $f$-orbit of $x$.}\\
In order to prove it, we note that for any $n\in \NN$, the relation $gf^n = f^n g$ implies
$$\left|\frac{\Jac^s(g)(f^{n} x)}{\Jac^s(g)(x)}\right| = 
\left|\frac{\Jac^s(f^{n})(g x)}{\Jac^s(f^{n})(x)}\right|.$$
Since $f^{n}(x)$ lies in a compact region of $W^s(\cO)$ for all $n\in \NN$,
the left hand side of this expression is uniformly bounded in $n$.
Hence if $y=g(x)$, then the quantity $\left|\log(\Jac^s(f^n)(x))-\log(\Jac^s(f^n)(y))\right|$
is bounded in $n$.
Since $g(\cO)=\cO$ the point $y=g(x)$ belongs to $W^s(\cO)\setminus \cO$ and
the stable manifold unbounded distortion property implies that $y$ belongs to the $f$-orbit of $x$,
proving the claim.

For $i\in \ZZ$ we consider the closed subset $P_i=\{x\in W^s(\cO)\setminus \cO, g(x)=f^i(x)\}$
of $W^s(\cO)\setminus \cO$.
The claim above shows that their union $P$ contains $\cD$, hence is dense in $W^s(\cO)\setminus \cO$.
Moreover for $i\neq j$ the intersection $P_i\cap P_j$ is empty
(any point in the intersection would be $j-i$-periodic).

We claim that:
\emph{For every $i\in \ZZ$, the set $P_i$ is open in $P$.}\\
Consider a point $x\in P_i$ and a neighborhood $U\subset W^s(\cO)$
of $x$ that is disjoint from all its iterates.
Consider also a point $y\in P\cap U$ close to $x$.
Since $g(x)$ belongs to $f^i(U)$, this is also the case for $g(y)$.
On the other hand there exists $j\in \ZZ$ such that $g(y)=f^j(y)$.
Since the iterates of $U$ are pairwise disjoint, it follows that $j=i$.
The map $g$ hence agrees with $f^i$ on a neighborhood of $x$ in $P\cap W^s(\cO)$,
proving the claim.

Let us consider a non-empty set $P_i$.
Since $P$ is dense in $W^s(\cO)$, it follows from the second claim that the interior of
$P_i$ in $W^s(\cO)$ is a non-empty $f$-invariant set in which $g=f^i$.
Since $\cO$ is a hyperbolic periodic orbit, the Lipschitz constant of $f^i$
in a neighborhood of $\cO$ is arbitrarily large when $|i|$ goes to infinity.
Using that $g$ is uniformly Lipschitz, one deduces that only a finite number of sets
$P_i$ are non-empty. Consequently $W^s(\cO)\setminus \cO$ is the disjoint 
union of finitely many closed sets $P_i$. Hence any connected component 
$W$ of $W^s(\cO)\setminus \cO$ is contained in a single set $P_i$, proving that $g=f^i$ on $W$.
\end{proof}

The periodic orbits of a generic diffeomorphism $f$ are fixed by $C(f)$.
\begin{lemma}\label{l.fixe}
For any diffeomorphism $f$ in a residual set $\cR_{per}$ of
$\Diff^1(M)$ (resp. of $\hbox{Symp}^1(M)$, of $\Diff^1_\mu(M)$),
for any periodic orbit $\cO$ of $f$ and for any $g\in \Lip(M)$
such that $f\circ g=g\circ f$, we have $g(\cO)=\cO$.
\end{lemma}
\begin{proof}
If $f$ and $g$ commute, then any periodic orbit $\cO$ of $f$
is sent by $g$ on a periodic orbit $\cO'$ of $f$ having the same period.
Since $g$ is Lipschitz, the eigenvalues of $\cO$ and $\cO'$ for $f$ must coincide.
The set $\cR_{per}$ of diffeomorphisms whose periodic orbits have different eigenvalues
hence satisfies the conclusion of the lemma.
\end{proof}

To prove Theorem~\ref{t=main}, it thus remains to show:
\begin{prop}\label{p=main}
For any diffeomorphism $f$ in a residual set $\cR\subset \Diff^1(M)$,
any periodic orbit of $f$ is hyperbolic and has the stable manifold unbounded distortion property.

For any diffeomorphism $f$ in a residual set $\cR_{symp}\subset \hbox{Symp}^1(M)$,
any hyperbolic periodic orbit of $f$ has the stable unbounded distortion property. 

For any diffeomorphism $f$ in a residual set $\cR_{\mu}\subset \Diff^1_\mu(M)$
any hyperbolic periodic $\cO$ of $f$ whose stable manifold $W^s(\cO)$ has codimension at least $\dim(M)/2$
has the stable manifold unbounded distortion property.
\end{prop}

\subsection{Reduction to contractions of $\RR^d$}
Let $B^d$ denote the unit closed ball $\overline{B(0,1)}$ of $\RR^d$
and consider the Banach space of $C^1$ maps $B^d\to\RR^d$
that send $0$ to $0$, endowed with the $C^1$ topology given by the $C^1$ norm:
$\|f-g\|_1=\sup_{x\in B^d} \|f(x)-g(x)\|+ \|D_xf-D_xg\|$.
The set of embeddings $B^d\to \RR^d$ fixing $0$
defines an open subset that will be denoted by $\cD^d$.

A \emph{contraction} of $\RR^d$ is an element of $\cD^d$ that sends $B^d$ into ${B(0,1)}$,
so that $0$ is a (hyperbolic) sink that attracts all the points in $B^d$.
The set of contractions of $\RR^d$ is an open subset   $\cC^d\subset
\cD^d$, and hence a Baire space.

Let $f$ be a diffeomorphism of a manifold $M$,  $p$ be a periodic point of $f$, $d^s$
be its stable dimension and $\tau$ be its period.
A \emph{stable chart} for $p$ is a local chart $\psi\colon \RR^d\to M$ such that,
denoting by $\pi^s$ the projection of $\RR^d$ onto the $d^s$ first coordinates,
we have the following properties.
\begin{itemize}
\item The domain $\psi(\RR^d)$ contains $p$.
\item In the chart $\psi$, the local stable manifold of $p$ contains the graph
of a $C^1$ map $g\colon \RR^{d^s}\to \RR^{d-d^s}$.
\item Let $v$ be equal to $\pi^s(\psi^{-1}(p))$ and let $\theta$ be the $C^1$ map
defined on a neighborhood of $0$ by projecting on the space $\RR^{d^s}$
the dynamics of $f^\tau$ in the local stable manifold of $p$:
$$\theta\colon x\mapsto \pi^s\circ \psi^{-1}\circ f^\tau \circ \psi(x+v, g(x+v)) -v.$$
Then, $\theta$ belongs to $\cC^{d^s}$.
\end{itemize}

The following property is easy to check.
\begin{lemma}\label{l.stable-chart}
Any hyperbolic periodic point $p$ of a diffeomorphism $f$ has a stable chart $\psi$.
Moreover, for any diffeomorphism $g$ in a $C^1$ neighborhood $\cU$ of $f$, the continuation
$p_g$ of $p$ also admits the chart $\psi$ as a stable chart.

The family of contractions $\theta_g$ associated to the periodic point $p_g$
and to the chart $\psi$ induces a continuous map $\Theta\colon \cU\to \cC^{d^s}$. 
This map is open.
\end{lemma}

In the conservative setting, the same property holds, but the proof is more delicate and
is postponed until section~\ref{s.invmanif}.
\begin{theorem}\label{t.stable-chart}
Let $\Theta\colon \cU\to \cC^{d^s}$ be a family of contractions associated to
a periodic point $p$ and a stable chart $\psi$ as in lemma~\ref{l.stable-chart}.
Then the map $\Theta\colon \cU\cap \hbox{Symp}^1(M)\to \cC^{d^s}$ is open.

If the dimension $d^s$ of the stable space of $p$ is greater than or equal to $\dim(M)/2$,
then the map $\Theta\colon \cU\cap \Diff_\mu^1(M)\to \cC^{d^s}$ is open.
\end{theorem}

The major ingredient in the proof of Proposition~\ref{p=main} is the following.

\begin{prop}\label{p=diskdist}
There is a residual set $\cS^{d} \subset \cC^d$
and a dense set $\cD\subset B^d\setminus \{0\}$ such that any $f\in \cS^d$
has the following property:
for any points $x\in \cD$, $y \in B^d\setminus \{0\}$ not in the same orbit
and for any $N\geq 1$ there exists $n\geq 1$ such that
\begin{eqnarray}\label{e=jacdist}
\left|\log(\Jac(f^n)(x))-\log(\Jac(f^n)(y))\right|>N.
\end{eqnarray}
\end{prop}

All these results together allow us to prove that $C^1$ generically, the stable manifold unbounded distortion
property holds.

\noindent
\begin{proof}[Proof of Proposition~\ref{p=main}]
For any integer $n\geq 0$, there exists
\begin{itemize}
\item a family $\cP_n$ of pairwise disjoint open subsets whose union is dense in $\Diff^1(M)$,
\item for each $\cU\in \cP_n$, finitely many charts $\psi_1,\dots,\psi_s\colon \RR^d\to M$,
\end{itemize}
such that any diffeomorphism $f\in \cU$ has the following properties:
\begin{itemize}
\item $f$ has $s$ periodic points of period less than $n$, all are hyperbolic.
Each domain $\psi_i(\RR^d)$ contains exactly one of them, it is called $p_{i,f}$
and its stable dimension is denoted by $d_i^s$ and its period by $\tau_i$.
\item The chart $\psi_i$ is a stable chart for $p_{i,f}$.
\end{itemize}
By Lemma~\ref{l.stable-chart}, the chart $\psi_i$ induces an open map $\Theta_i\colon \cU\to \cC^{d^s_i}$,
so that $\Theta_i^{-1}(\cS^{d^s_i})$ is residual in $\cU$ where $\cS^{d_i^s}$ is the residual subset of
$\cC^{d^s_i}$ provided by proposition~\ref{p=diskdist}.
Suppose that $f$ belongs to this residual set.
The dynamics of $f^{\tau_i}$ on the local stable manifold of $p_i$
are differentiably conjugate to the dynamics of a map $\Theta_i(f)\in \cS^{d^s_i}$,
proving that $p_i$ has the stable unbounded distortion property.

The union $\cR_n$ over $\cU\in \cP_n$
of the sets $\bigcap_{i=1}^s \Theta_i^{-1}(\cS^{d^s_i})$ is residual in $\Diff^1(M)$
and Proposition~\ref{p=main} holds with the residual set $\cR=\cap_n \cR_n$.
The proof in the conservative cases is similar and uses Theorem~\ref{t.stable-chart}.
\end{proof}

\subsection{Huge distortion: proof of Proposition~\ref{p=diskdist}}
We fix a countable dense subset $\cD$ of $\interior(B^d)\setminus \{0\}$ and a countable basis
of compact neighborhoods $\cB$ for $B^d\setminus \{0\}$.
For any $x\in \cD$ and any compact set $\Lambda\in \cB$,
the set of contractions $f\in \cC^d$ such that the orbits of $x$ and of $\Lambda$
are disjoint is an open set $O(x,\Lambda)$.
For any $N\geq 1$, the subset $O(x,\Lambda,N)\subset O(x,\Lambda)$ of contractions $f$
such that for any $y\in \Lambda$ we have~(\ref{e=jacdist}) is open.
Thus the set
$$\cS^d=\bigcap_{x,\Lambda,N} \left( O(x,\Lambda,N)\cup (\cC^d\setminus \overline{O(x,\Lambda)})\right)$$
is a $G_\delta$ and any diffeomorphism $f\in \cS^d$ satisfies the conclusion of Proposition~\ref{p=diskdist}.
In order to prove Proposition~\ref{p=diskdist}
we must show that $\cS^d$ is dense (hence residual) in $\cC^d$.
This is a consequence of the following perturbation result.

\begin{lemma}
The set $O(x,\Lambda,N)$ is dense in $O(x,\Lambda)$
\end{lemma}
\begin{proof}
Consider a contraction $f\in O(x,\Lambda)$.
By performing $C^1$ small perturbation, we can assume that $f$ is linear in a neighborhood $V$ of $0$.
Now consider $i_0\geq 0$ such that
all the iterates $f^i(x)$ and $f^i(\Lambda)$ with $i\geq i_0$ are contained in $V$.
There exists a $C^1$ small perturbation with compact support $g_1$ of the linear map $D_0f$
such that $\Jac g_1(0)\neq \Jac f(0)$.

For some $n\geq i_0$ and
any $i\in\{i_0,\dots n\}$, the map $f$ agrees with $D_0f$ in a neighborhood of $f^i(x)$
and can be replaced by the map
$$g_{i,\varepsilon}\colon z\mapsto f^i(x)+\varepsilon. g_1\left(\frac{z-f^i(x)}{\varepsilon}\right).$$
By choosing $\varepsilon$ small enough, all these perturbations have disjoint support.
We thus obtain a map $g$ that is still $C^1$ close to $f$.
Note that $f$ and $g$ agree on the positive orbit of $x$ and of a neighborhood of $K$.

We have for any point $y \in \Lambda$,
\begin{equation*}
\begin{split}
&\left|\log(\Jac(f^n)(x)) -\log(\Jac(f^n)(y))\right|\geq \\
&\quad\quad(n-i_0) \left|\log\Jac g_1(0)-\log \Jac f(0)\right|- \left|\log(\Jac(f^{i_0})(x))-\log(\Jac(f^{i_0})(y))\right|,
\end{split}
\end{equation*}
which is larger than $N$, if $n$ has been chosen large enough. Hence, $g$ belongs to
$O(x,\Lambda,N)$.
\end{proof}

\section{Conservative extension results}\label{s.invmanif}
We explain in this part how a perturbation of a conservative diffeomorphism
along a submanifold $W$ can be extended to a conservative perturbation on the whole manifold $M$.

This implies Theorem~\ref{t.stable-chart}:
the results proven in this section will be applied to the case where $W$ is an invariant manifold
of a hyperbolic periodic point $p$. In the volume-preserving case, we
assume that $\dim(W)\leq \frac 1 2 \dim(M)$ (note that this hypothesis is always satisfied either
by the stable or by the unstable manifold of $p$).
In the symplectic case, there is no additional hypothesis, but we use the following well-known fact.

\begin{lemma} Let $f\in \hbox{Symp}^1(M)$ and let $p$ be a hyperbolic periodic point for $f$.
Then $W^s(p)$ and $W^u(p)$ are Lagrangian submanifolds of $M$.
\end{lemma}
\begin{proof}
Let $x\in W^s(p)$, and let $v,w\in T_xW^s(p)$ be
tangent vectors to $W^s(p)$. On the one hand, since $f$ is a symplectomorphism,
we have
$$\omega(D_xf^k(v), D_xf^k(w)) = \omega(v,w),$$
for all $k\in \ZZ$.
On the other hand, as $k\to +\infty$, we have
$$\omega(D_xf^k(v), D_xf^k(w))\to 0.$$
Hence $\omega$ vanishes identically on $W^s(p)$.
The same is true for $W^u(p)$. Since $W^s(p)$ and $W^u(p)$ have
complementary dimension and $\omega$ is nondegenerate, they must have the same dimension.
Hence, both are Lagrangian submanifolds of $M$.
\end{proof}

\subsection{The symplectic case}\label{s=symplectic}
\begin{prop}\label{p.sympert} 
Let $M$ be a symplectic manifold and $z$ a point contained in a $C^1$ Lagrangian
submanifold $W\subset M$. Then there exists in $W$ a disk $D=\overline{B_W(z,r_0)}$ centered at $z$ such that,
for every neighborhood $U\subset M$ of $D$ and every $\eps>0$,
there exists $\delta>0$ with the following property.

For every $C^1$ diffeomorphism $\psi:D\to D$ satisfying:
\begin{itemize}
\item[a.] $\psi=Id$ on a neighborhood of $\partial D$, and  
\item[b.] $d_{C^1}(\psi, Id) <\delta$,
\end{itemize}
there exists $\varphi\in\hbox{Symp}^1(M)$ such that:
\begin{enumerate}
\item $\varphi=Id$ on $M\setminus U$,
\item $\varphi = \psi$ on $D$, and
\item $d_{C^1}(\varphi,Id) <\eps$.
\end{enumerate}
\end{prop}
\begin{proof}
The basic strategy is first to symplectically
embed the disk $D$ as the $0$-section of its cotangent bundle
$T^*D$. On $T^*D$, the symplectic form is $\omega = d\alpha$, where $\alpha$ is 
the canonical one-form on $T^*D$.  Any diffeomorphism $\psi:D\to D$ 
lifts to a canonical symplectomorphism $\psi^\ast: T^\ast D\to T^\ast D$;
namely the pull-back map $(\psi, D\psi^{-1})$.  The natural thing to
try to do is to set 
$\varphi = \psi^\ast$ in a neighborhood of the $0$-section, symplectically
interpolating between $\psi^\ast$ and $Id$ using a generating function.
This simple approach fails, however, because $\psi$ is only $C^1$,
and so $\psi^\ast$ is merely continuous.
(Even assuming that $\psi$ is $C^2$ does not help:
in order to control the $C^1$ size
of such a map, it is necessary to have some control on the $C^2$ size of
$\psi$, and we cannot assume any such control).  Using a convolution, 
it is possible to overcome this problem.   This approach mirrors that in 
\cite{BGV}, but in the symplectic setting.

The problem is local and one can work in $\RR^{2n}$ endowed with the standard symplectic form
$\omega = \sum_i du_i\wedge dv_i$ where $u = (u_1,\ldots, u_n), v=(v_1,\ldots, v_n)$.
By a symplectic change of coordinates, we may assume that the disk $D$ lies inside
a disk $\{(u,v),\,\|u\|\leq R, \, v=0\}$.
We define $\psi$ using a generating function $S$.  

We first recall the definition and properties of generating functions.
Suppose that $h:{\RR^{2n}}\to \RR^{2n}$ is a $C^r$ symplectomorphism, taking the form:
$$h(u,v) = (\xi(u,v),\eta(u,v)),$$
with $\xi,\eta:{\RR}^{2n}\to \RR^n$ and
$h(0,0) = (0,0)$.
Assume that the partial derivative matrix $\frac{\partial}{\partial v}\eta(u,v)$ is
invertible (this is the case for instance if $h$ preserves $\RR^n \times \{0\}$).
We can solve for $\eta = \eta(u,v)$ to obtain new coordinates
$(u,\eta)$ on a small neighborhood of $(0,0)$ in  $\RR^{2n}$.
Since $h$ is symplectic, the $1$-form $\alpha = \sum_i v_i du_i + \xi_i d\eta_i$ is closed, and hence, exact.
Thus there exists a $C^{r+1}$ function $S = S(u, \eta)$, unique up to adding a constant, defined in
a neighborhood of $(0,0)$, such that $dS = \alpha$.
The function $S$ is called a {\em generating function for $h$}.

On the other hand, any $C^{r+1}$ function $S=S(u,\eta)$ satisfying the nondegeneracy condition 
that $\frac{\partial^2}{\partial u\partial \eta}S$ is everywhere nonsingular is the generating 
function of a $C^r$ symplectic diffeomorphism. Solving for $\alpha$ in the
equation
$$dS = \frac{\partial S}{\partial u} du + \frac{\partial S}{\partial\eta} d\eta = \alpha = v du + \xi d\eta ,$$
we obtain the system:
$$\frac{\partial S}{\partial u} = v;\qquad \frac{\partial S}{\partial\eta} = \xi.$$
The nondegeneracy condition implies that this system can be solved
implicitly for a $C^r$ function $\eta = \eta(u,v)$.  
We then obtain a  $C^r$ symplectomorphism:
$$h(u,v) = \left(\frac{\partial S}{\partial\eta}(u,\eta(u,v)), \eta(u,v)\right),$$
and $S$ is a generating function for $h$.

It is easy to see that the generating function for the identity map is
$$S_0(u,\eta) = u\cdot \eta = \sum_{i=1}^n u_i\eta_i.$$

\begin{affi}
For every $\eps>0$, there exists $\delta>0$ such that, 
if $d_{C^2}(S,S_0) < \delta$ then $d_{C^1}(h,Id) <\eps$.
\end{affi}

\begin{proof} This follows from the implicit function theorem,
and the details are omitted.
\end{proof}

Returning to the proof of Proposition~\ref{p.sympert}, assume  
that $\psi:D\to D$ is written in $u$-coordinates as
$$\psi(u_1,\ldots, u_n) = (\psi_1(u_1,\ldots,u_n),\ldots, \psi_n(u_1,\ldots, u_n)).$$
We may assume that the domain of $\psi$ has been extended to 
$\RR^n$.  
To prove Proposition~\ref{p.sympert}, it suffices to prove the following lemma.
\end{proof}

\begin{lemma}\label{l.sympert} Given a disk $D\subset \RR^n$,
and a neighborhood $U$ of $D\times \{0\}$
in $\RR^{2n}$,  there exists $C>0$ with the following property.

For every $C^1$ diffeomorphism $\psi:\RR^n\to \RR^n$, equal
to the identity on a neighborhood of $\partial D$, there is a $C^2$
function $S:\RR^{2n}\to \RR$ such that:
\begin{enumerate}
\item $d_{C^2}(S_0,S) \leq Cd_{C^1}(\psi, Id)$,
\item $S=S_0$ outside of $U$,
\item $\frac{\partial S}{\partial u}(u,0) = 0$ for all $u\in \RR$ and
\item $\frac{\partial S}{\partial \eta}(u,0) = \psi(u)$ for all $u\in D$.
\end{enumerate}
\end{lemma}
Note that condition 1. implies that $S$ is nondegenerate, provided that
$d_{C^1}(\psi, Id)$ is sufficiently small.

\bigskip

\noindent{\bf Proof of Lemma~\ref{l.sympert}.}
To illustrate the argument in a simple case, we
first prove the lemma for $n=1$.
The proof of the general case is very similar.
Let
$$a(u) = \psi'(u)-1.$$ 
Note that $a$ is a continuous map, $\|a\|_\infty\leq d_{C^1}(\psi, Id)$, and
$a(u) =0$ if $u\notin \hbox{int}(D)$. 
Let $\Phi:\RR\to [0,1]$ be a $C^\infty$ function satisfying:
\begin{itemize}
\item $\Phi(0) = 1$ and $\Phi = 0$ outside of $(-1,1)$,
\item $\Phi^{(k)}(0) = 0$, for all $k\geq 1$,
\item $\int_\RR \Phi(w)\, dw = 1$.
\end{itemize}
Fix a point $u_\ast\in \partial D$, so that $\psi(u_\ast) = u_\ast$.
For $(u,\eta)\in \RR^{2}$, $\eta\neq 0$, let:
$$Q(u,\eta) = \eta \int_{u_\ast}^u\int_\RR \Phi(w)\, a(x-w\eta)\,dw\,dx.$$
For $\eta\neq 0$, one can make the change of variables
$w' = x-w\eta$ and get
$$Q(u,\eta)= \sign(\eta)\int_{u_\ast}^u \int_\RR
\Phi\left(\frac{x-w'}{\eta}\right) a(w')\,dw'\,dx.$$
Let $\rho:\RR^2\to \RR$ be a $C^\infty$ bump function
identically equal to $1$ on a neighborhood of $D\times \{0\}$
and vanishing outside of $U$. Consider
$$S= S_0 + \rho\, Q.$$
Lemma~\ref{l.sympert} in the case $n=1$ is a direct consequence of:

\begin{clai}\label{l.conv} The map $Q\colon \RR^{2}\to \RR$ is $C^2$ and 
there is $C = C(U) >0$ such that:
\begin{enumerate}
\item $\|Q\,\vert_{\overline U}\|_{C^2} \leq C \|a\|_\infty$,
\item $\frac{\partial Q}{\partial u}(u,0) = 0$, for all $u\in \RR$, and
\item $\frac{\partial Q}{\partial \eta}(u,0) = 
\int_{u_\ast}^u a(x)\,dx = \psi(u)-u$, for all  $u\in \RR$.
\end{enumerate}
\end{clai}
\begin{proof}
We derive explicitly the formulas: 
\begin{eqnarray*}
\frac{\partial Q}{\partial u} 
&=& \eta \int_\RR \Phi(w) \, a(u-w\eta)\,dw\\
&=& \sign(\eta)\int_{\RR} \Phi\left(\frac{u-w'}{\eta}\right)\, a(w')\,dw',\\
\\
\frac{\partial Q}{\partial \eta}& =&
\frac{-\sign(\eta)}{\eta^2}\int_{u_\ast}^u \int_{\RR} \Phi'\left(\frac{x-w'}{\eta}\right)\,(x-w') \,a(w')\,dw'\,dx\\
&=& -\int_{u_\ast}^u \int_\RR \Phi'(w)\, w\, a(x-w\eta)\, dw\,dx\\
&=&-\int_\RR \Phi'(w)\, w\,\int_{u_\ast - w\eta}^{u-w\eta} a(x')\,dx'\,dw,
\end{eqnarray*}
\begin{eqnarray*}
\frac{\partial^2 Q}{\partial \eta \partial u}& =& 
-\int_\RR \Phi'(w)\, w\, a(u-w\eta)\, dw,\\
\\
\frac{\partial^2 Q}{\partial u^2} &=&
\frac{1}{|\eta|}\int_{\RR} \Phi'\left(\frac{u-w'}{\eta}\right) \,a(w')\,dw'\\
&=& \int_\RR \Phi'(w)\,a(u -w\eta)\, dw,\\
\\
\hbox{and finally:}&&\\
\\
\frac{\partial^2 Q}{\partial \eta^2}
&=& \int_\RR \Phi'(w) \,w^2 \,\left(a(u-w\eta) - a(u_\ast -w\eta)\right) \,dw.
\end{eqnarray*}
Properties 1. and 2. follow immediately 
from these formulas. To see 3., note that
\begin{eqnarray*}
\frac{\partial Q}{\partial \eta}\vert_{\eta = 0}
&=& -\left(\int_{u_\ast}^u a(x)\,dx\right)  \left(\int_\RR \Phi'(w) w\, dw\right)\\
&=& -\left(\int_{u_\ast}^u a(x)\,dx\right)  \left(-\int_\RR \Phi(w) dw\right)\\
&=&\int_{u_\ast}^u a(x)\,dx.
\end{eqnarray*}
\end{proof}

We now turn to the case $n\geq 1$ in Lemma~\ref{l.sympert}.
For $i=1,\ldots n$, let $\alpha_i$ be the continuous $1$-form defined by
$$\alpha_i = d(\psi_i - \pi_i),$$
where $\pi_i:\RR^n\to \RR$ is the projection onto the $i$th coordinate.
As above, fix a point $u_\ast\in \partial D$, so that $\psi(u_\ast) = u_\ast$.
Then we have the formula:
$$\psi_i(u_1,\ldots,u_n) - u_i =
\int_{u_\ast}^u \alpha_i,$$
where the right-hand side is a path integral evaluated
on any path from $u_\ast$ to $u = (u_1,\ldots, u_n)$.
Furthermore, we have $\|\alpha_i\|_\infty \leq d_{C^1}(\psi,Id)$, for all $i$.
When $n=1$, the $1$-form $\alpha_1$ is just
$\alpha_1 = a(u)\, du$, where $a(u) = \psi'(u)-1$, as above.

Let $\Phi_n:\RR^n\to\RR$ be an 
$n$-dimensional bell function:
$$\Phi_n(x_1,\ldots, x_n) = \Phi(x_1)\cdots \Phi(x_n).$$
For each $1$-form $\alpha$, and $t\in \RR$, we define
a new $1$-form $\alpha_i^{\star t}$ on $\RR^n$ by taking the convolution:
$$\alpha^{\star t}(u) = t\int_{\RR^n}  \Phi_n(w)\, \alpha(u-tw) \, dw.$$

We integrate along any path from $u_\ast$ to $u$ and set
$$Q(u,\eta) = \int_{u_\ast}^u \sum_{i=1}^n \alpha_i^{\star \eta_i}
=\sum_{i=1}^n \eta_i\int_{\RR^n}\Phi_n(w) \,\left(\int_{u_\ast}^u
\alpha_i(u-tw)\right)\,dw.$$
This is well-defined since
$\int_{u_\ast}^u\alpha_i(u-tw)$
is independent of choice of path.

Let $\rho_n:{\bf R}^{2n} \to [0,1]$ be a $C^\infty$ bump
function vanishing identically outside of $U$ and equal to $1$
on a neighborhood of $D$.
As before, the map $S = S_0 +\rho_n Q$ satisfies the conclusions of Lemma~\ref{l.sympert}
provided the following claim holds.

\begin{clai} The map $Q\colon\RR^{2n}\to \RR$ is $C^2$ and there is $C = C(U) >0$ such that:
\begin{enumerate}
\item $\|Q\,\vert_{\overline U}\|_{C^2} \leq C \max_i\|\alpha_i\|_\infty$,
\item $\frac{\partial Q}{\partial u}(u,0) = 0$, for all $u\in \RR^n$, and
\item $\frac{\partial Q}{\partial \eta_i}(u,0) = 
\int_{u_\ast}^u \alpha_i  = \psi_i(u)-u_i$, for all $1\leq i\leq n$ and  $u\in \RR$.
\end{enumerate}
\end{clai}
\begin{proof}
We repeat the calculations from the proof of Lemma~\ref{l.conv} in the general setting.
When $t\neq 0$, the change of variable $w'=u-tw$ gives
\begin{eqnarray*}
\alpha^{\star t} (u) &=& \sign(t)\int_{\RR^n} \Phi_n
\left(u-tw'\right)\,\alpha(w' )\,dw',\\
\frac{d}{dt} \alpha^{\star t} (u) &=& -\int_{\RR^n}
\left(d\Phi_n(w).w+(n-1)\Phi_n(w)\right)\,\alpha(u-tw)\,dw.
\end{eqnarray*}
One deduces:
\begin{eqnarray*}
\frac{\partial Q}{\partial u} &=&
\sum_{i=1}^n \eta_i \int_{\RR^n} \Phi_n(w)\, \alpha_i(u-tw)\,dw\\
&=& \sum_{i=1}^n \frac{\eta_i}{|\eta_i|^n}\int_{\RR^n}
\Phi_n\left(\frac{u-w'}{\eta_i}\right)\, \alpha_i(w')\,dw',\\
\\
\frac{\partial Q}{\partial \eta_i}& =&
 \int_{u_\ast}^u \frac{d}{d\eta_i} \alpha_i^{\star\eta_i}\\
&=&-\int_{\RR^n} \left(d\Phi_n(w).w+(n-1)\Phi_n(w)\right)\,
\left(\int_{x=u_\ast-\eta_iw}^{x=u-\eta_iw}\alpha_i(w)\right) \, dw, \\
\\
\frac{\partial^2 Q}{\partial u\, \partial \eta_i}
&=&-\int_{\RR^n} \left(d\Phi_n(w).w+(n-1)\Phi_n(w)\right)\,
\alpha_i(u-\eta_iw)\, dw,\\
\\
\frac{\partial^2 Q}{\partial u^2} &=&
\sum_{i=1}^n 
\int_{\RR^n} d \Phi_n(w)\,\alpha_i(u-\eta_iw)\, dw,\\
\end{eqnarray*}
and finally:
\begin{eqnarray*}
\frac{\partial^2 Q}{\partial \eta_i \partial \eta_j}
&=&  \delta_{i,j} \int_{\RR^n} \left(d\Phi_n(w).w+(n-1)\Phi_n(w)\right)\,
\left[\alpha_i(x-\eta_iw).w\right]_{x=u_\ast}^{x=u} \, dw.\\
\end{eqnarray*}
It is not difficult to verify that 1.--3. hold.
\end{proof}

The proof of Lemma~\ref{l.sympert} is now complete.\endproof

\subsection{The volume-preserving case}\label{s=volpres}

\begin{prop}\label{p.vpert} Let $M$ be a Riemannian manifold
endowed with a volume form $\mu$ and $W$ be a $C^1$ submanifold satisfying
$$\hbox{dim}(W) \leq \hbox{codim}(W).$$
Centered at any point $z\in W$, there exists a disk $D=\overline{B_W(z,r_0)}$
of $W$ such that, for every neighborhood $U\subset M$ containing $D$ and every $\eps>0$,
there exists $\delta>0$ with the following property.  

For every $C^1$ diffeomorphism $\psi:D\to D$ satisfying:
\begin{itemize}
\item[a.] $\psi=Id$ on a neighborhood of $\partial D$, and  
\item[b.] $d_{C^1}(\psi, Id) <\delta$,
\end{itemize}
there exists $\varphi\in\Diff^1_\mu(M)$ such that:
\begin{enumerate}
\item $\varphi=Id$ on $M\setminus U$,
\item $\varphi = \psi$ on $D$, and 
\item $d_{C^1}(\varphi,Id) <\eps$.
\end{enumerate}
\end{prop}
\begin{proof} Let $n=\hbox{dim}(M)$.
By a local change of coordinates, we may assume that
$\mu$ is the standard volume form $dx_1\wedge \cdots\wedge dx_n$
on a neighborhood of the origin in $\RR^n$.
By composing these coordinates with an isometry of $\RR^n$,
we may further assume that $D$ is the graph of a $C^1$ function $h:\RR^k \to \RR^{n-k}$,
where $k\leq n/2$.
The final change of coordinates
$$(x_1,\ldots, x_n)\mapsto ((x_1,\ldots, x_k), (x_{k+1},\ldots, x_n) - h(x_1,\ldots, x_k))$$
preserves volume.  Applying this change of coordinates, we may assume
that $D$ lies in the coordinate plane
$\{x_{k+1}=x_{k+2}=\cdots=x_n=0\}\simeq \RR^k$. Now we apply the symplectic
pertubation result (Proposition~\ref{p.sympert}) inside the
space $\{x_{2k+1}=\cdots=x_n = 0\}\simeq \RR^{2k}$ to obtain a local
$C^1$ symplectomorphism $\varphi_0$ of $\{x_{2k+1}=\cdots=x_n= 0\}$ that agrees with 
$\psi$ on $D$. This symplectomorphism is 
$C^1$ isotopic to the
identity through symplectomorphisms $\{\varphi_t\}_{t\in[0,1]}$,
where $\varphi_1 = Id$ (to obtain this
isotopy, just choose a smooth
isotopy of the generating function for $\psi$
to the generating function for the identity).
 
Now we extend $\varphi_0$ to $\RR^n$
using this isotopy to obtain a 
locally-supported volume-preserving diffeomorphism
that agrees with $\psi$ on $D$.  More precisely, choose an
 appropriate
$C^\infty$ bump function $\rho:\RR^{n-2k}\to [0,1]$, and set
$$\varphi(x_1,\ldots, x_n) = (\varphi_{\rho(\|(x_{2k+1},\ldots,x_n)\|)}(x_1,\ldots, x_{2k}), x_{2k+1},\ldots, x_n).$$
This is the desired map $\varphi$.
\end{proof}

\vspace{10pt}

\noindent \textbf{Christian Bonatti (bonatti@u-bourgogne.fr)}\\
\noindent  CNRS - Institut de Math\'ematiques de Bourgogne, UMR 5584\\
\noindent  BP 47 870\\
\noindent  21078 Dijon Cedex, France\\
\vspace{10pt}

\noindent \textbf{Sylvain Crovisier (crovisie@math.univ-paris13.fr)}\\
\noindent CNRS - Laboratoire Analyse, G\'eom\'etrie et Applications, UMR 7539,\\
\noindent Institut Galil\'ee, Universit\'e Paris 13, 99 Avenue J.-B. Cl\'ement,\\
\noindent 93430 Villetaneuse, France\\
\vspace{10pt}

\noindent \textbf{Amie Wilkinson (wilkinso@math.northwestern.edu)}\\
\noindent Department of Mathematics, Northwestern University\\
\noindent 2033 Sheridan Road \\
\noindent Evanston, IL 60208-2730,  USA

\end{document}